\documentclass[a4paper,12pt]{article}
\usepackage{amssymb}
\usepackage{amsthm}
\usepackage{amsmath}
\usepackage[utf8]{inputenc}
\usepackage{array}
\usepackage[T1]{fontenc}
\usepackage{lmodern}
\usepackage{tikz}
\usepackage{vmargin}
\setmarginsrb{2.5cm}{1cm}{2.5cm}{1cm}{1cm}{1cm}{1cm}{1cm}

\title{Fractional triangle decompositions in graphs with large minimum degree}

\usepackage{hyperref}
\usepackage[affil-it]{authblk}

\author[a]{François Dross}

\affil[a]{{\small ENS de Lyon, LIRMM}}  \affil[ ]{{\small 161 rue Ada, 34095 Montpellier
    Cedex 5, France}} \affil[
]{\href{mailto:francois.dross@ens-lyon.fr}{\small{francois.dross@ens-lyon.fr}}}
\date{}

\begin{document}
\maketitle

\newtheorem{theo}{Theorem}
\newtheorem*{theo*}{Theorem}
\newtheorem{cor}[theo]{Corollary}
\newtheorem{lemm}[theo]{Lemma}
\newtheorem{prop}[theo]{Property}
\newtheorem{obs}[theo]{Observation}
\newtheorem{conj}[theo]{Conjecture}
\newtheorem{claim}[theo]{Claim}
\newtheorem{config}[theo]{Configuration}
\newtheorem{quest}[theo]{Question}

\begin{abstract}
A triangle decomposition of a graph is a partition of its edges into triangles. A fractional triangle decomposition of a graph is an assignment of a  non-negative weight to each of its triangles such that the sum of the weights of the triangles containing any given edge is one. We prove that for all $\epsilon > 0$, every large enough graph graph on $n$ vertices with minimum degree at least $(0.9 + \epsilon)n$ has a fractional triangle decomposition. This improves a result of Garaschuk that the same result holds for graphs with minimum degree at least $0.956n$. Together with a recent result of Barber, K\"{u}hn, Lo and Osthus, this implies that for all $\epsilon > 0$, every large enough triangle divisible graph on $n$ vertices with minimum degree at least $(0.9 + \epsilon)n$ admits a triangle decomposition.
\end{abstract}

\section{Introduction}

Decomposition and packing problems are central and classical problems in combinatorics, in particular, in design theory.
Kirkman's theorem \cite{kirkman} from the middle of 19th century gives a necessary and sufficient condition on the existence of a Steiner
triple system with a certain number of elements. In the language of graph theory, Kirkman's result asserts that
every complete graph with an odd number of vertices and a number of edges divisible by three can be decomposed into triangles.
Barber, K\"{u}hn, Lo and Osthus \cite{barber} showed that the same conclusion is true for large graphs satisfying necessary divisibility conditions if their minimum
degree is not too far from the number of their vertices. In this short paper, we study the fractional variant of the problem and
we use it to improve the bound obtained by Barber et al.

Let us fix the terminology we are going to use. A \emph{graph} is a pair of sets $\left(V,E\right)$ such that elements of $E$ are unordered pairs of elements of $V$. The elements of $V$ are called \emph{vertices} and the elements of $E$ are called \emph{edges}. We denote by $uv$ (or $vu$) the edge with vertices $u$ and $v$. We denote by $|G|$ the number of vertices of $G$. Two vertices contained in the same edge are said to be \emph{adjacent} or to be \emph{neighbours}. Two edges that share a vertex are said to be \emph{adjacent}. The \emph{degree} of a vertex $v$ is equal to the number of neighbours of $v$. Let $\gcd\left(G\right)$ denote the greatest common divisor of the degrees of the vertices of $G$.

Two graphs $G_1 = \left(V_1,E_1\right)$ and $G_2 = \left(V_2,E_2\right)$ are \emph{isomorphic} if there exists a bijection $b$ from $V_1$ to $V_2$ such that $uv$ is an edge of $G_1$ if and only if $b\left(u\right)b\left(v\right)$ is an edge of $G_2$ for every two vertices $u$ and $v$ of $G_1$.
The \emph{complete graph} $K_k$ is the graph with $k$ vertices all mutually adjacent. The graph $K_3$ is also called a \emph{triangle}. A graph $G_1 = \left(V_1,E_1\right)$ is a \emph{subgraph} of $G_2 = \left(V_2,E_2\right)$ if $V_1 \subseteq V_2$ and $E_1 \subseteq E_2$. The subgraphs of $G_2$ isomorphic to $G_1$ will be referred to as copies of $G_1$.

Let $H$ be a graph. An \emph{$H$-decomposition} of a graph $G$ is a set of subgraphs of $G$ isomorphic to $H$ that are edge disjoint such that each edge of $G$ is contained in one of them. A graph is \emph{$H$-decomposable} if it admits an $H$-decomposition. A $K_3$-decomposition is also called a \emph{triangle decomposition} and a graph is \emph{triangle decomposable} if it is $K_3$-decomposable. A graph $G$ is \emph{$H$-divisible} if $\gcd\left(G\right)$ is a multiple of $\gcd\left(H\right)$, and the number of edges of $G$ is a multiple of the number of edges of $H$. It is easy to see that every $H$-decomposable graph is $H$-divisible, but the converse is not true. Kirkman~\cite{kirkman} proved that every $K_3$-divisible complete graph is $K_3$-decomposable. The fact that for all $H$, every $H$-divisible complete graph is $H$-decomposable remained an open problem for over one hundred years before it was solved by Wilson \cite{wilson}.

The first generalisation to graphs that are near complete is due to Gustavsson \cite{gustavsson}. He proved that for every graph $H$, there exist $n_0(H)$ and $\epsilon(H)$ such that every $H$-divisible graph with $n \ge n_0(H)$ vertices and minimum degree at least $(1-\epsilon(H))n$ is $H$-decomposable. This has been generalised to hypergraphs in a recent result of Keevash \cite{keevash}. The best that is known to date for a general graph $H$ is due to Barber et al. \cite{barber}, who proved that for all $\epsilon > 0$, every sufficiently large $H$-divisible graph on $n$ vertices with minimum degree at least $\left(1-\frac{1}{16|H|^2(|H|-1)^2} + \epsilon\right)n$ is $H$-decomposable. For some particular classes of graphs, the exact asymptotic minimum degree threshold is known \cite{barber}\cite{yuster}.

A \emph{fractional $H$-decomposition} of a graph $G$ is an assignment of non-negative weights to the copies of $H$ in $G$ such that for an edge $e$, the sum of the weights of the copies of $H$ that contain $e$ is equal to one. A graph is \emph{fractionally $H$-decomposable} if it admits a fractional $H$-decomposition. A graph can be fractionally $H$-decomposable without being $H$-divisible. A fractional $K_3$-decomposition is also called a \emph{fractional triangle decomposition} and a graph is \emph{fractionally triangle decomposable} if it is fractionally $K_3$-decomposable. For all $r \ge 2$, Yuster \cite{yusterbis} proved that every graph on $n$ vertices with minimum degree at least $\left(1 - \frac{1}{9r^{10}}\right)n$ is fractionally $K_r$-decomposable, and Dukes \cite{dukes} proved that the same result holds for sufficiently large graphs on $n$ vertices with minimum degree at least $\left(1-\frac{1}{16r^2(r-1)^2}\right)n$.

In this paper we will focus on triangle decompositions of graphs with large minimum degree.
The following conjecture is due to Nash-Williams~\cite{nashwilliams}:

\begin{conj}[Nash-Williams~\cite{nashwilliams}] \label{c1}
Let $G$ be a $K_3$-divisible graph with $n$ vertices and minimum degree at least $\frac{3}{4}n$. If $n$ is large enough, then $G$ is $K_3$-decomposable.
\end{conj}

The best result towards a proof of Conjecture~\ref{c1} is due to Barber et al.~\cite{barber}.

\begin{theo}[Barber et al.~\cite{barber}] \label{bar}
There exists an $n_0$ such that every $K_3$-divisible graph $G$ on $n \ge n_0$ vertices with minimum degree at least $0.956n$ is $K_3$-decomposable.
\end{theo}

The proof of Theorem~\ref{bar} relies on a result on fractional $K_3$-decomposability, which we now state. The following appears as a conjecture in \cite{garaschuk}.

\begin{conj}[Garaschuk~\cite{garaschuk}] \label{c2}
Let $G$ be a graph with $n$ vertices and minimum degree at least $\frac{3}{4}n$. If $n$ is large enough, then $G$ is fractionally $K_3$-decomposable.
\end{conj}

The best known result towards proving Conjecture~\ref{c1} was established by Garaschuk \cite{garaschuk}.

\begin{theo}[Garaschuk~\cite{garaschuk}] \label{gar}
Let $G$ be a graph with $n$ vertices and minimum degree at least $0.956n$. The graph $G$ admits a fractional triangle decomposition.
\end{theo}

In this paper we use a different method to prove the following.

\begin{theo} \label{main}
Let $\epsilon > 0$. There exists an $n_0$ such that every graph with $n \ge n_0$ vertices and minimum degree at least $(\frac{9}{10} + \epsilon)n$ admits a fractional triangle decomposition.
\end{theo}

In \cite{barber}, a particular case of Theorem 11.1 and Lemma 12.3 imply the following.

\begin{theo}[Barber et al.~\cite{barber}] \label{barbis}
Suppose there exist $n_0$ and $\delta$ such that every graph on $n \ge n_0$ vertices with minimum degree at least $\delta n$ is fractionally $K_3$-decomposable. For all $\epsilon > 0$, there exist $n_1$ such that every $K_3$-divisible graph on $n \ge n_1$ vertices with minimum degree at least $\left(\max\left(\delta,\frac{3}{4}\right) + \epsilon\right)n$ vertices is $K_3$-decomposable.
\end{theo}

Together with Theorem \ref{barbis}, our result improves Theorem~\ref{bar}.

\begin{theo} \label{cmain} 
Let $\epsilon > 0$. There exists an $n_0$ such that every $K_3$-divisible graph on $n \ge n_0$ vertices with minimum degree at least $(\frac{9}{10} + \epsilon)n$ is $K_3$-decomposable.
\end{theo}

\section{Proof of Theorem~\ref{main}} \label{proofmain}

Let $\delta<\frac{1}{10}$ and fix a graph $G$ with $n$ vertices and minimum degree at least $\left(1-\delta\right)n$. Suppose the graph $G$ has at least one triangle with three vertices of degree at least $\left(1-\delta\right)n + 2$. Let $G'$ be the graph $G$ where the edges of one such triangle are removed.
Observe that $G'$ has minimum degree at least $\left(1-\delta\right)n$ and that if $G'$ has a fractional triangle decomposition, then $G$ has one too. Up to doing this operation several times, we can assume that $G$ has no triangle with three vertices of degree at least $\left(1-\delta\right)n + 2$.
Let $m$ be the number of edges of $G$. 

Initially, we give the same weight $w_\Delta$ to every triangle such that the sum of the weights of the triangles is equal to $3m$. We will modify the weights of the triangles to obtain a fractional triangle decomposition. We will do so in a way that the total sum of the weights is preserved. 

\begin{figure}[t]
\begin{center}
\begin{tikzpicture}
\coordinate (v0) at (0,0) ;
\coordinate (v1) at (0,2) ;
\coordinate (v2) at (1.772,-1) ;
\coordinate (v3) at (-1.772, -1) ; 
\coordinate (w1) at (0,-1) ;
\coordinate (w2) at (0.866,0.5) ;
\coordinate (w3) at (-0.866, 0.5) ; 
\coordinate (u1) at (0,1) ;
\coordinate (u2) at (0.866,-0.5) ;
\coordinate (u3) at (-0.866, -0.5) ; 

\draw[very thick] (v0) -- (v1);
\draw (v0) -- (v2);
\draw (v0) -- (v3);
\draw (v1) -- (v2);
\draw (v1) -- (v3);
\draw[dashed] (v2) -- (v3);

\draw [fill=black] (v0) circle (1.5pt) ;
\draw [fill=black] (v1) circle (1.5pt) ; 
\draw [fill=black] (v2) circle (1.5pt) ; 
\draw [fill=black] (v3) circle (1.5pt) ;

\draw (w1) node [below] {$-$} ;
\draw (w1) node [above] {$-$} ;

\draw (w2) node [left] {$+$} ;
\draw (w2) node [right] {$-$} ;

\draw (w3) node [right] {$+$} ;
\draw (w3) node [left] {$-$} ;

\draw (u1) node [left] {$+$} ;
\draw (u1) node [right] {$+$} ;

\draw (u2) node [above] {$+$} ;
\draw (u2) node [below] {$-$} ;

\draw (u3) node [above] {$+$} ;
\draw (u3) node [below] {$-$} ;

\end{tikzpicture}
\caption{By removing some weight $w$ from two triangles and adding $w$ to the two other triangles, we remove $2w$ from the dashed edge and add $2w$ to the thick edge. \label{ftrick}}
\end{center}
\end{figure}
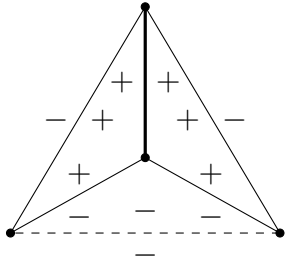

We define the weight of an edge $e$ to be the sum of the weights of the triangles that contain $e$. Given $H$ a copy of $K_4$ in $G$, and two non-adjacent edges $e_1$ and $e_2$ in $H$, let us call $\left(H, \{e_1,e_2\}\right)$ a \emph{rooted $K_4$} of $G$. We will use the following procedure to modify the weights of the edges of a rooted $K_4$ of $G$:

\emph{
Let $\left(H, \{e_1,e_2\}\right)$ be a rooted $K_4$ of $G$. By removing a weight $w$ from the two triangles of $H$ that contain $e_1$ and adding the same weight $w$ to each of the other two triangles (i.e. those that contain $e_2$), we transfer a weight of $2w$ from $e_1$ to $e_2$. The weights of all the other edges of the graph remain unchanged (see Figure~\ref{ftrick}).
}

To prevent the weight of any triangle from becoming negative, we have to restrict how much weight we can transfer using the procedure above. If for some $w$ we use the procedure to transmit a weight of $2w$ from an edge to another one, then any triangle's weight is lowered by at most $w$ for triangles that are in the $K_4$, and does not change for other triangles.
Moreover, since every triangle contains a vertex with degree at most $(\delta + 1)n$, any triangle is in at most $(1-\delta)n$ copies of $K_4$, and thus in at most $(1-\delta)n$ oriented copies of $K_4$ (since for each $K_4$ there are three possible choices for the pair of edges). Since each triangle has an initial weight of $w_\Delta$, if it sends weights of at most $\frac{2w_\Delta}{3(1-\delta)n}$ through each rooted $K_4$ that it is contained in, its final weight will be non-negative. 

We express redistributing the weights as a flow problem in an auxiliary graph, which is denoted by $\widehat G$. The vertices of $\widehat G$ are the edges of $G$. Two vertices in $\widehat G$ are adjacent if they form a pair in a rooted $K_4$ and the edge between them is set to have the capacity $c = \frac{2w_\Delta}{3(1-\delta)n}$. Let $E_c$ be the set of these edges. Moreover $\widehat G$ has two additional vertices, which we will call the {\em supersource} and the {\em supersink}. Let $T_e$ be the number of triangles of $G$ that contains an edge $e$. If $T_ew_\Delta>1$, then the vertex of $\widehat G$ corresponding to $e$ is joined to the supersource and the capacity of the corresponding edge of $\widehat G$ is $T_ew_\Delta - 1$. Likewise, if $T_ew_\Delta<1$, then the vertex of $\widehat G$ corresponding to $e$ is joined to the supersink and the capacity of the corresponding edge is $1 - T_ew_\Delta$. The vertices of $G$ adjacent to the supersource are referred to as {\em sources} and those adjacent to the supersink as {\em sinks}. Let
$$M = \sum_{e \text{ source}} \left(T_ew_\Delta - 1\right) = \sum_{e \text{ sink}} \left(1 - T_ew_\Delta\right).$$
We will show that $\widehat G$ has a flow of value $M$ from the supersource to the supersink.

If $\widehat G$ does not have a flow of value $M$, then it has a vertex cut $\left(A_0,B_0\right)$ such that the supersource is contained in $A_0$, the supersink in $B_0$ and the sum of the capacities of the edges from $A_0$ to $B_0$ is less than $M$. Let $A$ be the edges of $G$ corresponding to the vertices of $A_0$ and $B$ the edges corresponding to the vertices of $B_0$. Note that $|A|=|A_0|-1$ and $|B|=|B_0|-1$. Finally, let $k=|A|$, and observe that $|B|=m-k$.

Let $T_A$ and $T_B$ be the average $T_e$ for $e$ in $A$ and in $B$ respectively. Let $e = uv$ be an edge of $G$. Let $W_e$ be the set of the vertices $w$ such that $uvw$ is a triangle. By the definition of $T_e$, $|W_e| = T_e$. Each vertex of $W_e$ is non-adjacent to at most $\delta n$ vertices of $G$, and thus is non-adjacent to at most $\delta n$ vertices of $W_e$. So each vertex of $W_e$ is adjacent to at least $T_e - \delta n$ vertices of $W_e$. Therefore $e$ is in at least $\frac{T_e\left(T_e - \delta n\right)}{2}$ distinct copies of $K_4$, and consequently $e$ is in at least $\frac{T_e\left(T_e - \delta n\right)}{2}$ rooted $K_4$.

Let $e$ be a vertex of $A$. It is adjacent to at least $\frac{T_e\left(T_e - \delta n\right)}{2} - k$ vertices of $B$. Therefore the cut contains at least $$\sum_{e\in A}\left(\frac{T_e\left(T_e - \delta n\right)}{2} - k\right)$$ edges of $E_c$. Similarly, it contains at least $$\sum_{e\in B}\left(\frac{T_e\left(T_e - \delta n\right)}{2} - \left(m-k\right)\right)$$ edges of $E_c$. Moreover, for each source $e$ that is in $B$ and each sink $e$ that is in $A$, the cut contains the edge between $e$ and the supersource or the supersink. Recall that the capacities of the edges of $E_1$ is $c = \frac{2w_\Delta}{3(1-\delta)n}$. Therefore the sum of the capacities of the edges of $\widehat G$ is at least $$\sum_{e\in A}\left(\frac{T_e\left(T_e - \delta n\right)}{2} - k\right)c + \sum_{e \text{ source} \in B}\left(T_ew_\Delta - 1\right) + \sum_{e \text{ sink} \in A}\left(1 - T_ew_\Delta\right).$$ At the same time, it is also at least
$$\sum_{e\in B}\left(\frac{T_e\left(T_e - \delta n\right)}{2} - \left(m-k\right)\right)c + \sum_{e \text{ source} \in B}\left(T_ew_\Delta - 1\right) + \sum_{e \text{ sink} \in A}\left(1 - T_ew_\Delta\right).$$ Since the sum of the capacities of the edges in the considered cut is less than $M$, we get that
\begin{equation}
\sum_{e\in A}\left(\frac{T_e\left(T_e - \delta n\right)}{2} - k\right)c + \sum_{e \text{ source} \in B}\left(T_ew_\Delta - 1\right) + \sum_{e \text{ sink} \in A}\left(1 - T_ew_\Delta\right) < M \label{eq1}
\end{equation}
and
\begin{equation}
\sum_{e\in B}\left(\frac{T_e\left(T_e - \delta n\right)}{2} - \left(m-k\right)\right)c + \sum_{e \text{ source} \in B}\left(T_ew_\Delta - 1\right) + \sum_{e \text{ sink} \in A}\left(1 - T_ew_\Delta\right) < M. \label{eq2}
\end{equation}

The inequalities (\ref{eq1}) and (\ref{eq2}) can be rewritten using that $$M = \sum_{e \text{ source}} \left(T_ew_\Delta - 1\right)$$ and $$M = \sum_{e \text{ sink}} \left(1 - T_ew_\Delta\right)$$ respectively as follows.

\begin{equation}\sum_{e\in A}\left(T_e\left(T_e - \delta n\right) - 2k\right)c - 2\sum_{e \in A} \left(T_ew_\Delta - 1\right) < 0 \label{neq1} \end{equation}

\begin{equation}\sum_{e\in B}\left(T_e\left(T_e - \delta n\right) - 2(m-k)\right)c - 2\sum_{e \in B} \left(1 - T_ew_\Delta\right) < 0 \label{neq2} \end{equation}

Since the summand is a convex function of $T_e$, we obtain the following.

\begin{equation}
(T_A(T_A-\delta n)-2k)c-2(T_Aw_\Delta-1)<0 \label{feq1}
\end{equation}

\begin{equation}
(T_B(T_B-\delta n)-2(m-k))c-2(1-T_Bw_\Delta)<0 \label{feq2}
\end{equation}

The inequality (\ref{feq1}) implies that
\begin{equation}
T_A(T_A-\delta n) + \frac{2}{c}(1-T_Aw_\Delta) < 2k. \label{feq1bis}
\end{equation}

The inequality (\ref{feq2}) implies that
\begin{equation}
2k < 2m - T_B(T_B-\delta n) + \frac{2}{c}(1-T_Bw_\Delta). \label{feq2bis}
\end{equation}

We now combine the inequalities (\ref{feq1bis}) and (\ref{feq2bis}) and we substitute $c = \frac{2w_\Delta}{3(1-\delta)n}$ to get the following
\begin{equation}
T_A(T_A-\delta n) -(3n(1-\delta)T_A) < 2m - T_B(T_B-\delta n) -(3n(1-\delta)T_B) \label{eq3}
\end{equation}

Let $e$ be an edge of $G$. Each end-vertex of $e$ is non-adjacent to at most $\delta n$ vertices of $G$. Hence, the edge $e$ is contained in at least $n-2\delta n$ triangles. Since $e$ cannot be contained in more than $n$ triangles, we get that $n-2\delta n \le T_e \le n$. Consequently, we have $n-2\delta n \le T_A,T_B \le n$.

A standard analytic argument shows that the left hand side of (\ref{eq3}) is minimized when $T_A = n$ and the right hand side is maximized when $T_B = n - 2 \delta n$. Consequently, it must hold that 
\begin{equation}
n(n-\delta n) -3(1-\delta)n^2 < 2m - (n-2\delta n)(n-3\delta n) -(3n(1-\delta)(n-2\delta n)) \label{neq3}
\end{equation}

Note that we proved that in $G$, there is no triangle with three vertices of degree at least $(1-\delta)n+2$. Since $G$ has minimum degree at least $(1-\delta)n$, this implies that there are at most $2\delta n$ vertices of degree at least $(1-\delta)n+2$. Therefore we have $2m \le (2\delta n)n + ((1 - 2\delta)n)((1-\delta)n+1)$

We get that $1 - 11\delta + 10\delta^2 < \frac{1- 2\delta}{n}$. Since $\delta < \frac{1}{10}$, if $n$ is large enough this leads to a contradiction. Since there exists a flow of value $M$ in $\widehat G$, the weights of the triangles can be redistributed in a way that the triangles form a fractional decomposition of $G$. This finishes the proof of Theorem \ref{main}.

\section{Conclusion}

In this paper we proved that for all $\epsilon > 0$, there exists a constant $n_0$ such that every graph on $n \ge n_0$ vertices with minimum degree at least $\left(\frac{9}{10} + \epsilon\right)n$ is fractionally triangle decomposable. This implies that for all $\epsilon > 0$, there exists a constant $n_0$ such that every triangle divisible graph on $n \ge n_0$ vertices with minimum degree at least $\left(\frac{9}{10} + \epsilon\right)n$ is triangle decomposable.

\section{Acknowledgements}
This work was done during the author's visit to the group of Dan Kr\'al' at the University of Warwick; the visit was partially supported from the European Research Council under the European Union's Seventh Framework Programme (FP7/2007-2013)/ERC grant agreement no.~259385. I am deeply grateful to Tereza Klimo\v{s}ov\'{a} for her helpful comments, to Dan Kr\'al' for his careful reading and many suggestions, and to Ben Barber and Richard Montgomery for pointing out a significant improvement.

\bibliographystyle{plain}
\bibliography{biblio} {}

\end{document}